\documentclass[a4paper,oneside,10pt,,mathscr,final]{amsart}

\usepackage{amssymb}
\usepackage{euscript}
\usepackage[frame,cmtip,curve,arrow,matrix,line,graph]{xy}

\newcommand{\bb} {\mathbb}

\newcommand{\cal}{\mathcal}
\newcommand{\frk}{\mathfrak}

\newcommand{\wtld}{\widetilde}
\newcommand{\wh}{\widehat}

\newcommand{\ep}{\epsilon}

\newcommand{\e}{\mathrm{e}}
\newcommand{\tw}{\mathrm{tw}}
\newcommand{\Id}{\mathrm{Id}}

\newcommand{\Aut}{\operatorname{Aut}}
\newcommand{\Auteq}{\operatorname{Auteq}}
\newcommand{\Coh}{\operatorname{Coh}}
\newcommand{\Ext}{\operatorname{Ext}}
\newcommand{\Ho} {\operatorname{Ho}}
\newcommand{\Hom}{\operatorname{Hom}}
\newcommand{\Iso}{\operatorname{Iso}}

\newcommand{\Obj}{\operatorname{Obj}}
\newcommand{\Pic}{\operatorname{Pic}}
\newcommand{\SL}{\operatorname{SL}}

\newcommand{\blob}{{\scriptscriptstyle \bullet}}
\newcommand{\cl}{\widehat}
\newcommand{\longinto}{\lhook\joinrel\longrightarrow}

\newcommand{\<}{\langle}
\renewcommand{\>}{\rangle}

\newcommand{\bu}[1]{#1_{\blob}}
\newcommand{\cpx}[4]{
 \begin{xy}
  \xymatrix@=0pt@C=20pt@R=15pt{#1
  \ar@<.5ex>[r]^{#2} & \ar@<.5ex>[l]^{#3}  #4}
 \end{xy}  
}

\theoremstyle{plain}
 \newtheorem{thm}{Theorem}[section]

 \newtheorem{cor}[thm]{Corollary}
 \newtheorem*{thm*}{Theorem}
\theoremstyle{definition}
 \newtheorem{dfn}[thm]{Definition}
 \newtheorem{fct}[thm]{Fact}
 \newtheorem{df}[thm]{Definition/Fact}
\theoremstyle{remark}
 \newtheorem{rmk}[thm]{Remark}

\numberwithin{equation}{section}
\allowdisplaybreaks

\begin{document}


\title{A note on Bridgeland's Hall algebra of two-periodic complexes}
\author{Shintarou Yanagida}
\address{Research Institute for Mathematical Sciences,
Kyoto University, Kyoto 606-8502, Japan}
\email{yanagida@kurims.kyoto-u.ac.jp}

\thanks{The author is supported by JSPS Fellowships 
for Young Scientists (No.\ 21-2241, 24-4759)}

\date{July 4, 2012}


\begin{abstract}
We show that the Hall algebra of two-periodic complexes,
which is recently introduced by T.~Bridgeland,
coincides with the Drinfeld double 
of the ordinary Hall bialgebra.
\end{abstract}

\maketitle


\section{Introduction}


\subsection{}

The main object of this paper is the Hall algebra of 
$\bb{Z}_2(:=\bb{Z}/2\bb{Z})$-graded complexes,
which was introduced by Bridgeland \cite{B}.
 
Let $\cal{A}$ be an abelian category over a finite field $\frk{k} := \bb{F}_q$
with finite dimensional morphism spaces.
Let $\cal{P} \subset \cal{A}$ be the subcategory of projective objects.
Let $\cal{C}(\cal{A}) \equiv \cal{C}_{\bb{Z}_2}(\cal{A})$ be 
the abelian category of $\bb{Z}_2$-graded complexes in $\cal{A}$.
An object of $\cal{C}(\cal{A})$ is of the form
\begin{align*}
\cpx{M_1}{f}{g}{M_2},
\quad
f\circ g =0,\ 
g\circ f =0.
\end{align*}
Let $\cal{C}(\cal{P})$ be the subcategory of complexes 
consisting of projectives,
and $\cal{H}(\cal{C}(\cal{P}))$ be its Hall algebra.
One can introduce the twisted Hall algebra 
$\cal{H}_{\tw}(\cal{C}(\cal{P}))$ 
as the twisting of $\cal{H}(\cal{C}(\cal{P}))$ by the Euler form of $\cal{A}$. 
In \cite{B}, 
Bridgeland introduced an algebra $\cal{DH}(\cal{A})$,
which is the localization of the twisted Hall algebra  
$\cal{H}_{\tw}(\cal{C}(\cal{P}))$ by the set of acyclic complexes:
\begin{align*}
\cal{DH}(\cal{A}) := 
\cal{H}_{\tw}(\cal{C}(\cal{P}))
\bigl[ [\bu{M}]^{-1} \mid H_*(\bu{M})=0\bigr].
\end{align*}

The purpose of this note is to show the following theorem,
which was stated in \cite[Theorem~1.2]{B}.

\begin{thm*}
Assume that the abelian category $\cal{A}$ satisfies 
the conditions 
\begin{itemize}
\item essentially small with finite morphism spaces,
\item linear over $\frk{k}$,
\item of finite global dimension and having enough projectives,
\item hereditary,
\item nonzero object defines nonzero class in the Grothendieck group.
\end{itemize}
Then the algebra $\cal{DH}(\cal{A})$ 
is isomorphic to the Drinfeld double of 
the bialgebra $\wtld{\cal{H}}(\cal{A})$
as an associative algebra.
\end{thm*}

Here $\wtld{\cal{H}}(\cal{A})$ 
is the (ordinary) extended Hall bialgebra, 
which will be recalled in \S\S\ref{ssubsec:dfn} -- \ref{ssubsec:ba-hp}.
For the review of the Drinfeld double, see \S\ref{ssubsec:dd}.
The proof will be explained in \S\ref{sec:proof}.
\\

The organization of this note is as follows.
In \S\,\ref{sec:sec1}, we review Bridgeland's theory \cite{B} 
and prepare notations and statements 
which are necessary for the proof of the main theorem.
In the subsection \S\ref{subsec:Hall},
we recall the ordinary theory of Hall algebra 
introduced by Ringel \cite{R}.
The next subsection \S\ref{subsec:hac} is devoted 
to the recollection of Bridgeland's theory.

The section \S\ref{sec:proof} is devoted to the proof of the main theorem.

We close this note by mentioning some consequences of the theorem 
in \S\ref{sec:cr}.


\subsection{Notations and conventions}
\label{subsect:notation}

We indicate several global notations.

$\frk{k} := \bb{F}_q$ is a fixed finite field 
unless otherwise stated,
and all the categories will be $\frk{k}$-linear.
We choose and fix a square root $t:=\sqrt{q}$.

For an abelian category $\cal{A}$, 
we denote by $\Obj(\cal{A})$ the class of objects of $\cal{A}$.
For an object $M$ of $\cal{A}$, 
the class of $M$ in the Grothendieck group $K(\cal{A})$ 
is denoted by $\cl{M}$.
Let $\cal{K}_{\ge 0}(\cal{A}) \subset K(\cal{A})$ 
be the subset of $K(\cal{A})$ consisting 
of the classes $\cl{A} \in K(\cal{A})$ of $A \in \cal{A}$
(rather than the formal differences of them).

For an abelian category $\cal{A}$ 
which is essentially small, 
the set of its isomorphism classes is 
denoted by $\Iso(\cal{A})$.

For a complex 
$\bu{M}=(\cdots \to M_{i}\xrightarrow{d_i} M_{i+1} \to \cdots )$ 
in an abelian category $\cal{A}$,
its homology is denoted by $H_*(\bu{M})$. 

For a set $S$, 
we denote by $|S|$ its cardinality.


\section{Hall algebras of complexes}
\label{sec:sec1}


\subsection{Hall algebra}
\label{subsec:Hall}

This subsection gives basic definitions and properties of Hall algebra,
following \cite[\S\S2.3--2.5]{B} and \cite[\S1]{S}.
Let $\cal{A}$ be an abelian category satisfying the assumptions 
\begin{enumerate}
\item[(a)] essentially small with finite morphism spaces,
\item[(b)] linear over $\frk{k}$,
\item[(c')] of finite global dimension.
\end{enumerate}

\begin{rmk}
We will introduce additional conditions 
(c), (d) and (e) in the following discussion.
These symbols except (c') follow those in \cite{B}.
\end{rmk}


\subsubsection{Definitions}
\label{ssubsec:dfn}

Consider a vector space
$$
 \cal{H}(\cal{A}) := 
 \bigoplus_{A \in \Iso(\cal{A})} \bb{C}[A]
$$
linearly spanned by symbols $[A]$ with $A$ running through 
the set $\Iso(\cal{A})$  
of isomorphism classes of objects in $\cal{A}$.

\begin{df}[{Ringel \cite{R}}]

The following operation defines on $\cal{H}(\cal{A})$ 
the structure of a unital associative algebra over $\bb{C}$:
\begin{align}
\label{eq:sec1:diam}
[A] \diamond [B] := 
\sum_{C \in \Iso(\cal{A})} 
\dfrac{ |\Ext^1_{\cal{A}}(A,B)_{C}| }{ |\Hom_{\cal{A}}(A,B)|} \cdot [C].
\end{align}
Here 
\begin{align*}
\Ext^1_{\cal{A}}(A,B)_C \subset \Ext^1_{\cal{A}}(A,B)
\end{align*}
is the set parametrizing extensions of $B$ by $A$
with the middle term isomorphic to $C$.

The unit is given by $[0]$,
where $0$ is the zero object of $\cal{A}$.

This algebra $(\cal{H}(\cal{A}), \diamond, [0])$ is called 
the Hall algebra of $\cal{A}$. 
Below we will denote it by $\cal{H}(\cal{A})$ for simplicity.
\end{df}

\begin{rmk}
We follow \cite{B} to choose 
$|\Ext^1_{\cal{A}}(A,B)_{C}| / |\Hom_{\cal{A}}(A,B)| $
for the structure constant of the multiplication. 
It is proportional to the usual structure constant 
$|\{ B' \subset C \mid B' \cong B, C/B' \cong A  \}|$
appearing in \cite{R} and \cite{S}.
See \cite[\S2.3]{B} for the detail.
\end{rmk}


\subsubsection{Euler form and extended Hall algebra}
\label{ssubsec:eha}

Let us recall the notations for Grothendieck group 
given in \S\ref{subsect:notation}.
For objects $A,B\in \cal{A}$,
the Euler form is defined by 
\begin{align}\label{eq:sec1:euler}
\<A,B\> := \sum_{i \in \bb{Z}} (-1)^i \dim_{\frk{k}}\Ext^i_{\cal{A}}(A,B),
\end{align}
where the sum is finite by our assumptions on $\cal{A}$.
As is well known, this form descends to the one 
on the Grothendieck group $K(\cal{A})$ of $\cal{A}$,
which is denoted by the same symbol as \eqref{eq:sec1:euler}:
\begin{align*}
\<\cdot,\cdot\>: 
 K(\cal{A})\times K(\cal{A}) \longrightarrow \bb{Z}.
\end{align*}
We will also use the symmetrized Euler form:
\begin{align*}
(\cdot,\cdot): 
 K(\cal{A})\times K(\cal{A}) \longrightarrow \bb{Z},
 \quad
(\alpha,\beta) 
 := \<\alpha,\beta\> + \<\beta,\alpha\>.
\end{align*}

\begin{df}
\begin{enumerate}
\item[(1)]
The twisted Hall algebra $\cal{H}_{\tw}(\cal{A})$ 
is the same vector space as $\cal{H}(\cal{A})$  
with the twisted multiplication
\begin{equation*}
[A]*[B] := t^{ \< \cl{A},\cl{B} \> }\cdot [A] \diamond [B] 
\end{equation*}
for $A,B \in \Iso(\cal{A})$.
Here $t:=\sqrt{q}$ is the fixed square root of $q$.

\item[(2)]
The extended Hall algebra $\wtld{\cal{H}}(\cal{A})$ 
is defined as an extension of $\cal{H}_{\tw}(\cal{A})$
by adjoining symbols $K_\alpha$ for classes $\alpha\in K(\cal{A})$, 
and imposing relations
\begin{align*}
K_\alpha * K_\beta = K_{\alpha+\beta}, \quad 
K_\alpha * [B]     = t^{( \alpha,\cl{B})} \cdot [B] * K_\alpha
\end{align*}
for $\alpha,\beta \in K(\cal{A})$ 
and $B \in \Iso(\cal{A})$.
Note that $\wtld{\cal{H}}(\cal{A})$ 
has a vector space basis consisting of the elements 
$K_\alpha * [B]$ for $\alpha\in K(\cal{A})$ and $B \in \Iso(\cal{A})$. 
\end{enumerate}
\end{df}

\begin{rmk}
In \cite{B},
the extended Hall algebra is denoted by $\cal{H}_{\tw}^{\e}(\cal{A})$.
\end{rmk}


\subsubsection{Green's coproduct}
\label{ssubsec:cpr}

To introduce a coalgebra structure,
one should consider a completion of the algebra.
Assume that the abelian category $\cal{A}$ satisfies 
the conditions (a), (b), (c') and 
\begin{enumerate}
\item[(e)] nonzero object 
defines nonzero class in the Grothendieck group.
\end{enumerate}
Then the algebra $\cal{H}(\cal{A})$ is naturally graded 
by the Grothendieck group $K(\cal{A})$ of $\cal{A}$:
\begin{align*}
 \cal{H}(\cal{A}) 
 = \bigoplus_{\alpha \in K(\cal{A})}
 \cal{H}(\cal{A})[\alpha],
 \quad
 \cal{H}(\cal{A})[\alpha]:=\bigoplus_{\cl{A}=\alpha}\bb{C} [A].
\end{align*}
For $\alpha,\beta \in K(\cal{A})$, 
set
\begin{align*}
&\cal{H}(\cal{A})[\alpha] 
  \wh{\otimes}_{\bb{C}} 
 \cal{H}(\cal{A})[\beta]
 :=\prod_{\cl{A}=\alpha, \, \cl{B}=\beta} 
    \bb{C} [A] \otimes_{\bb{C}} \bb{C} [B],
 \\
&\cal{H}(\cal{A}) \wh{\otimes}_{\bb{C}} \cal{H}(\cal{A})
 :=\prod_{\alpha,\beta \in K(\cal{A})}
   \cal{H}(\cal{A})[\alpha] 
    \wh{\otimes}_{\bb{C}} \cal{H}(\cal{A})[\beta].
\end{align*}
Thus $\cal{H}(\cal{A}) \wh{\otimes} \cal{H}(\cal{A})$
is the space of all formal linear combinations
$$
 \sum_{A,B} c_{A,B}  \cdot [A] \otimes [B].
$$
This tensor product $\wh{\otimes}$ is called a completed tensor product.

\begin{df}
\begin{enumerate}
\item[(1)] (Green \cite{G})
The following maps
$$
 \Delta:\cal{H}(\cal{A}) \longrightarrow 
 \cal{H}(\cal{A}) \wh{\otimes}_{\bb{C}} \cal{H}(\cal{A}), \quad 
 \ep:\cal{H}(\cal{A}) \longrightarrow \bb{C}
$$ 
define a topological coassociative coalgebra structure
on $\cal{H}(\cal{A})$:
\begin{align}
\label{eq:Green}
 \Delta([A]) :=
  \sum_{B,C} t^{\<B,C\>} 
  \dfrac{|\Ext_{\cal{A}}(B,C)_A|}{|\Aut_{\cal{A}}(A)|}  \cdot 
   [B] \otimes [C],
 \quad
 \ep([A]):=\delta_{A,0}.
\end{align} 

\item[(2)] (Xiao \cite{X})
On the extended algebra, 
defining the maps 
$$
\Delta:\wtld{\cal{H}}(\cal{A}) \longrightarrow 
 \wtld{\cal{H}}(\cal{A}) \wh{\otimes}_{\bb{C}} \wtld{\cal{H}}(\cal{A}),
 \quad
 \ep: \wtld{\cal{H}}(\cal{A}) \longrightarrow \bb{C}
$$
by \eqref{eq:Green} and
\begin{align*}
 \Delta(K_\alpha) :=
  K_\alpha \otimes K_\alpha,
 \quad
 \ep(K_\alpha):=1,
\end{align*} 
one has a topological coassociative coalgebra structure
on $\wtld{\cal{H}}(\cal{A})$.
\end{enumerate}
\end{df}

Here the word \emph{topological} means that 
everything should be considered 
in the completed space.
For example, the coassociativity in (1) means that 
the two maps $(\Delta \otimes 1) \otimes \Delta$ 
and $(1 \otimes \Delta) \otimes \Delta$ 
from $\cal{H}(\cal{A})$ to 
$\cal{H}(\cal{A})\wh{\otimes}\cal{H}(\cal{A})
 \wh{\otimes}\cal{H}(\cal{A})$ 
coincide.

\subsubsection{Bialgebra structure and Hopf pairing}
\label{ssubsec:ba-hp}

Now we have an algebra structure 
and a coalgebra structure on $\cal{H}(A)$ 
(and on $\wtld{\cal{H}}(A)$).
In order that these structures are compatible and 
give a \emph{bialgebra} structure,
we must impose one more condition on $\cal{A}$.

\begin{fct}[{Green \cite{G}, Xiao \cite{X}}]
Assume that the abelian category $\cal{A}$ satisfies 
the conditions (a), (b), (c'), (e) and 
\begin{enumerate}
\item[(d)]
hereditary, that is, of global dimension at most 1.
\end{enumerate}
Then the tuples
$$
(\cal{H}(\cal{A}),\diamond,[0],\Delta,\ep),
\quad
(\wtld{\cal{H}}(\cal{A}),*,[0],\Delta,\ep)
$$
are topological bialgebras defined over $\bb{C}$.
That is,the map 
$\Delta:\cal{H}(\cal{A}) \to 
 \cal{H}(\cal{A}) \wh{\otimes} \cal{H}(\cal{A})$
and 
 $\Delta:\wtld{\cal{H}}(\cal{A}) \to 
 \wtld{\cal{H}}(\cal{A}) \wh{\otimes} 
 \wtld{\cal{H}}(\cal{A})$
are  homomorphisms of $\bb{C}$-algebras.
\end{fct}

Below, we simply denote by 
$\cal{H}(\cal{A})$ and $\wtld{\cal{H}}(\cal{A})$
the bialgebras
$(\cal{H}(\cal{A}),\diamond,[0],\Delta,\ep)$
and 
$(\wtld{\cal{H}}(\cal{A}),*,[0],\Delta,\ep)$
respectively.
\\

This bialgebra structure on $\cal{H}(\cal{A})$
has an additional feature,
that is, it is \emph{self-dual}.
The self-duality is stated in terms of 
a natural nondegenerate bilinear form, 
called \emph{Hopf pairing}.

\begin{df}[{Green, \cite{G}}]
\label{df:sec1:hp}
Assume that the abelian cateogry $\cal{A}$ satisfies 
the conditions (a), (b), (c'), (d) and (e).
\begin{enumerate}
\item[(1)]
The non-degenerate bilinear form 
$$
 (\cdot,\cdot)_H: 
 \cal{H}(\cal{A}) \otimes_{\bb{C}} \cal{H}(\cal{A})
 \longrightarrow
 \bb{C}
$$
given by 
$$
 ([A],[B])_H := \dfrac{\delta_{A,B}}{|\Aut_{\cal{A}}(A)|}
$$
is a Hopf pairing on the bialgebra $\cal{H}(\cal{A})$,
that is, 
for any $x,y,z \in \Iso(\cal{A})$,
one has 
\begin{equation}\label{eq:sec1:prdpair}
   (x \diamond y, z)_H
 = (x \otimes y, z)_H.
\end{equation}

\item[(2)]
The non-degenerate bilinear form 
$$
 (\cdot,\cdot)_H: 
 \wtld{\cal{H}}(\cal{A}) \otimes_{\bb{C}} 
 \wtld{\cal{H}}(\cal{A})
 \longrightarrow
 \bb{C}
$$
given by 
$$
 ([A] K_\alpha,[B] K_\beta)_H 
  := \dfrac{\delta_{A,B}}{|\Aut_{\cal{A}}(A)|}t^{(\alpha,\beta)}
$$
is a Hopf pairing 
on the bialgebra $\wtld{\cal{H}}(\cal{A})$,
that is, 
for any $x,y,z \in \Iso(\cal{A})$,
one has 
\begin{equation}\label{eq:sec1:prdpair2}
   (x * y, z)_H
 = (x \otimes y, z)_H.
\end{equation}
\end{enumerate}
\end{df}

\begin{rmk}
In the right hand sides of 
\eqref{eq:sec1:prdpair} and \eqref{eq:sec1:prdpair2},
we used the usual pairing on the product space:
$$
 (x \otimes y, z \otimes w)_H := 
 (x,z)_H \cdot (y, w)_H
$$
\end{rmk}

\subsubsection{Drinfeld double}
\label{ssubsec:dd}

Here we recall the Drinfeld double 
of the self-dual bialgebra.
For the complete treatment of 
Drinfeld double construction,
we refer \cite[\S3.2]{J} and \cite[\S5.2]{S}.

\begin{fct}[{Drinfeld}]
\label{fct:sec1:dd}
Let $\cal{H}$ be a $\bb{C}$-bialgebra 
with a Hopf pairing 
$(\cdot,\cdot)_H:\cal{H} \otimes_{\bb{C}} \cal{H} \to \bb{C}$.
Then there is a unique algebra structure $\circ$ on 
$\cal{H} \otimes \cal{H}$ 
satisfying the following conditions
\begin{enumerate}
\item[(1)]
The maps 
$$
 \cal{H} \longrightarrow \cal{H} \otimes_{\bb{C}} \cal{H},
 \quad
 a \longmapsto a\otimes 1
$$
and
$$
 \cal{H} \longrightarrow \cal{H} \otimes_{\bb{C}} \cal{H},
 \quad
 a \longmapsto 1 \otimes a
$$
are injective homomorphisms of $\bb{C}$-algebras.
\item[(2)]
For all elements $a,b \in \cal{H}$, one has
$$
 (a\otimes 1)\circ(1\otimes b)=a \otimes b.
$$
\item[(3)]
For all elements $a,b \in \cal{H}$, one has
\begin{equation}
\label{eq:sec1:dd}
  \sum (a_{(2)},b_{(2)})_H \cdot a_{(1)} \otimes b_{(1)}
 =\sum (a_{(1)},b_{(1)})_H \cdot
    (1 \otimes b_{(2)}) \circ (a_{(2)} \otimes 1).
\end{equation}
Here we used Sweedler's notation:
$\Delta(a)=\sum a_{(1)} \otimes a_{(2)}$.
\end{enumerate}
\end{fct}

\begin{rmk}
If $\cal{H}$ is a topological bialgebra,
then one should replace the tensor product $\otimes$ 
in the statement by the completed one $\wh{\otimes}$.
\end{rmk}


\subsection{Hall algebras of complexes}
\label{subsec:hac}
We summarize necessary definitions and properties of  
Hall algebras of $\bb{Z}_2$-graded complexes.
Most of the materials were introduced or shown in \cite{B}.

In this subsection \S\ref{subsec:hac}, 
$\cal{A}$ denotes an abelian category satisfying 
the following three conditions.
\begin{enumerate}
\item[(a)] essentially small with finite morphism spaces,
\item[(b)] linear over $\frk{k}$,
\item[(c)] of finite global dimension and having enough projectives.
\end{enumerate}

\subsubsection{Categories of two-periodic complexes}

We shall recall the basic definitions in \cite[\S3.1]{B}.
Let $\cal{C}_{\bb{Z}_2}(\cal{A})$ be the abelian category 
of $\bb{Z}_2$-graded complexes in $\cal{A}$.  
An object $\bu{M}$ of this category consists of 
the following diagram in $\cal{A}$:
\begin{align*}
\cpx{M_1}{d_1}{d_0}{M_0}, \quad d_{i+1}\circ d_i=0.
\end{align*} 
Hereafter indices in the diagram of an object 
in $\cal{C}_{\bb{Z}_2}(\cal{A})$ 
are understood by modulo $2$.
A morphism $\bu{s}: \bu{M} \to \bu{N}$ 
consists of a diagram
\begin{align*}
\xymatrix{
  M_1 \ar[d]_{s_1}\ar@<.5ex>[r]^{d_1} 
&     \ar@<.5ex>[l]^{d_0} M_0\ar[d]^{s_0}
\\
  N_1\ar@<.5ex>[r]^{d'_1}
&    \ar@<.5ex>[l]^{d'_0} N_0}
\end{align*}
with $s_{i+1} \circ d_{i} = d'_i\circ s_{i}$.
Two morphisms $\bu{s}, \bu{t}: \bu{M} \to \bu{N}$ 
are said to be homotopic 
if there are morphisms $h_i: M_i\to N_{i+1}$ such that
\begin{align*}
t_i-s_i=d'_{i+1} \circ h_{i}+h_{i+1}\circ d_{i}.
\end{align*}
For an object $\bu{M} \in C_{\bb{Z}_2}(\cal{A})$,
we define its class in the $K$-group by
$$
 \bu{\cl{M}} := \cl{M}_0 - \cl{M}_1 \in K(\cal{A}).
$$
\\

Denote by $\Ho_{\bb{Z}_2}(\cal{A})$ the category 
obtained from $\cal{C}_{\bb{Z}_2}(\cal{A})$ 
by identifying homotopic morphisms.
Let us also denote by 
\begin{align*}
\cal{C}_{\bb{Z}_2}(\cal{P}) \subset \cal{C}_{\bb{Z}_2}(\cal{A}),
\end{align*}
the full subcategories whose objects are complexes of projectives in $\cal{A}$.
Hereafter we drop the subscript $\bb{Z}_2$
and just write 
$$
 \cal{C}(\cal{A}) := \cal{C}_{\bb{Z}_2}(\cal{A}), \quad
 \cal{C}(\cal{P}) := \cal{C}_{\bb{Z}_2}(\cal{P}), \quad 
 \Ho(\cal{A})     := \Ho_{\bb{Z}_2}(\cal{A}).
$$
\\

The shift functor $[1]$ of complexes induces an involution 
\begin{align*}
\cal{C}(\cal{A}) \stackrel{*}{\longleftrightarrow} 
\cal{C}(\cal{A}).
\end{align*}
This involution shifts the grading and changes the sign of the differential 
as follows:
\begin{align*}
\bu{M} = \Bigl(\cpx{M_1}{d_1}{d_0}{M_0}\Bigr)
\ \stackrel{*}{\longleftrightarrow} \ 
\bu{M}^{*} = \Bigl(\cpx{M_0}{-d_0}{-d_1}{M_1}\Bigr)
\end{align*}

Now let us recall 

\begin{fct}[{\cite[Lemma~3.3]{B}}]\label{fct:B:lem3.3}
For $\bu{M}, \bu{N} \in \cal{C}(\cal{P})$ we have
\begin{align*}
\Ext_{\cal{C}(\cal{A})}^1(\bu{N},\bu{M}) \cong 
\Hom_{\Ho(\cal{A})}(\bu{N},\bu{M}^*).
\end{align*}
\end{fct}


A complex $\bu{M} \in \cal{C}(\cal{A})$ is called acyclic if $H_*(\bu{M})=0$.
To each object $P\in\cal{P}$, we can attach acyclic complexes
\begin{align*}
\bu{{K_P}}   = \bigl(\cpx{P}{\Id}{0}{P}\bigr), \qquad 
\bu{{K_P}}^* = \bigl(\cpx{P}{0}{-\Id}{P}\bigr).
\end{align*}

\begin{rmk}
The complexes $\bu{{K_P}}, \bu{{K_P}}^*$ 
are denoted by $K_P$, $K_P^*$ in \cite{B}.
\end{rmk}

Let us recall the following fact shown in \cite{B}.

\begin{fct}[{\cite[Lemma~3.2]{B}}]
For each acyclic complex of projectives $\bu{M} \in \cal{C}(\cal{P})$, 
there are objects $P,Q\in \cal{P}$, unique up to isomorphism, 
such that $\bu{M} \cong \bu{{K_P}} \oplus \bu{{K_Q}}^*$. 
\end{fct}


\subsubsection{Definition of Hall algebras of complexes}

Let $\cal{H}(\cal{C}(\cal{A}))$ be the Hall algebra of the abelian category $\cal{C}(\cal{A})$
defined in \S\ref{subsec:Hall}.
As noted in \cite[\S3.5]{B}, this definition makes sense 
since the spaces $\Ext^1_{\cal{C}(\cal{A})}(\bu{N},\bu{M})$ 
are all finite-dimensional by Fact~\ref{fct:B:lem3.3}.
Let
\begin{align*}
\cal{H}(\cal{C}(\cal{P}))\subset \cal{H}(\cal{C}(\cal{A}))
\end{align*}
be the subspace spanned by complexes of projective objects.

Define $\cal{H}_{\tw}(\cal{C}(\cal{P}))$ 
to be the same vector space as $\cal{H}(\cal{C}(\cal{P}))$ 
with the twisted multiplication
\begin{equation}
\label{eq:sec1:*}
[\bu{M}]*[\bu{N}]
:=t^{\<M_0,N_0\>+\<M_1,N_1\>} \cdot [\bu{M}] \diamond [\bu{N}].
\end{equation}

Now let us recall the simple relations satisfied 
by the acyclic complexes $\bu{{K_P}}$:

\begin{fct}[{\cite[Lemmas~3.4, 3.5]{B}}]
\label{fct:B:lem3.4-5}
For any object $P\in \cal{P}$ and any complex $\bu{M} \in \cal{C}(\cal{P})$, 
we have the following relations in $\cal{H}_{\tw}(\cal{C}(\cal{P}))$:
\begin{align*}
&[\bu{{K_P}}]*[\bu{M}] 
 = t^{ \<\cl{P},\bu{\cl{M}}\>} \cdot [\bu{{K_P}} \oplus \bu{M}],
\quad
 [\bu{M}]*[\bu{{K_P}}] 
 = t^{-\<\bu{\cl{M}},\cl{P}\>} \cdot [\bu{{K_P}} \oplus \bu{M}],
\\
&[\bu{{K_P}}]*[\bu{M}] 
 = t^{(\cl{P},\bu{\cl{M}})} \cdot [\bu{M}] * [\bu{{K_P}}],
\quad
 [\bu{{K_P}}^*]*[\bu{M}]
 =t^{-(\cl{P},\bu{\cl{M}})} \cdot [\bu{M}] * [\bu{{K_P}}^*]
\end{align*}
In particular, for $P,Q \in \cal{P}$ we have 
\begin{align}
\label{eq:additive}
\begin{split}
&[\bu{{K_P}}] * [\bu{{K_Q}}] 
   = [\bu{{K_P}} \oplus \bu{{K_Q}}], \quad 
 [\bu{{K_P}}] * [\bu{{K_Q}}^*] 
   = [\bu{{K_P}} \oplus \bu{{K_Q^*}}],
\\
&[[\bu{{K_P}}],[\bu{{K_Q}}]] 
   = [[\bu{{K_P}}],[\bu{{K_Q}}^*]]
   = [[\bu{{K_P}}^*],[\bu{{K_Q}}^*]]=0. 
\end{split}
\end{align}
At the last line  we used the commutator $[x,y] := x*y-y*x$.
\end{fct}

\subsubsection{Bridgeland's Hall algebra}

Now we can introduce the main object:
Bridgeland's Hall algebra.
We define the localized Hall algebra $\cal{DH}(\cal{A})$ 
to be the localization of $\cal{H}_{\tw}(\cal{C}(\cal{P}))$ 
with respect to the elements $[\bu{M}]$ 
corresponding to acyclic complexes $\bu{M}$: 
\begin{align*}
 \cal{DH}(\cal{A})
 := \cal{H}_{\tw}(\cal{C}(\cal{P}))
    \bigl[ [\bu{M}]^{-1} \mid H_*(\bu{M})=0\bigr].
\end{align*}
As explained in \cite[\S3.6]{B},
this is the same as localizing by the elements 
$[\bu{{K_P}}]$ and $[\bu{{K_P^*}}]$ 
for all objects $P\in \cal{P}$. 

For an element $\alpha \in K(\cal{A})$, we define 
\begin{align*}
K_\alpha := [\bu{{K_P}}] * [\bu{{K_Q}}]^{-1},\quad
K_\alpha^* := [\bu{{K_P}}^*] * [\bu{{K_Q}}^*]^{-1},
\end{align*}
where we expressed $\alpha = \cl{P}-\cl{Q}$ 
using the classes of some projectives $P,Q \in \cal{P}$.
This is well defined by Fact~\ref{fct:B:lem3.4-5}

\begin{rmk}
We will denote two different elements 
$K_\alpha \in \wtld{\cal{H}}(\cal{A})$ 
and 
$K_\alpha \in \cal{DH}(\cal{A})$ 
by the same symbol,
following \cite{B}.
\end{rmk}

By Fact \ref{fct:B:lem3.4-5}, we immediately have

\begin{cor}
\label{cor:sec1:kkc}
In the algebra $\cal{DH}(\cal{A})$,
we have 
\begin{enumerate}
\item[(1)]
$$
 K_{\alpha} * \bu{M} 
 = t^{(\alpha,\bu{\wh{M}})} \cdot \bu{M} * K_{\alpha},
 \quad
 K_{\alpha}^* * \bu{M} 
 = t^{-(\alpha,\bu{\wh{M}})} \cdot \bu{M} * K_{\alpha},
$$
for arbitrary $\alpha\in K(\cal{A})$
and $\bu{M} \in \cal{C}(\cal{P})$.

\item[(2)]
$$
 [K_\alpha,K_\beta] 
 = [K_\alpha,K_\beta^*]
 = [K_\alpha^*,K_\beta^*]=0
$$
for arbitrary $\alpha,\beta \in K(\cal{A})$.
\end{enumerate}
\end{cor}


\subsection{Hereditary case}
\label{subsec:sec1:hered}

In this subsection \S\ref{subsec:sec1:hered}, 
we assume that 
$\cal{A}$ satisfies the conditions (a), (b), (c) and 
the following additional ones:
\begin{itemize}
\item[(d)] $\cal{A}$ is hereditary, that is of global dimension at most 1,
\item[(e)] nonzero objects in $\cal{A}$ define nonzero classes in $K(\cal{A})$.
\end{itemize}
Then by \cite[\S4]{B} we have a nice basis for $\cal{DH}(\cal{A})$.
To explain that, let us recall the minimal resolution of objects of $\cal{A}$.

\subsubsection{Minimal resolution and the complex $\bu{{C_A}}$}

\begin{dfn}[{\cite[\S4.1]{B}}]
Assume the conditions (a),(c),(d) on $\cal{A}$. 
Then every object $A\in \cal{A}$ has a projective resolution
\begin{align}\label{eq:res}
0 \to P \xrightarrow{f} Q \xrightarrow{g} A \to 0,
\end{align}
and decomposing $P$ and $Q$ into finite direct sums
$P=\oplus_i P_i$, $Q=\oplus_j Q_j$,
one may write $f=(f_{ij})$ in matrix form with $f_{ij}: P_i\to Q_j$. 
The resolution \eqref{eq:res} is said to be minimal 
if none of the morphisms $f_{ij}$ is an isomorphism. 
\end{dfn}

\begin{fct}[{\cite[Lemma\,4.1]{B}}]\label{fct:B:lem4.1}
Any resolution \eqref{eq:res} is isomorphic to a resolution of the form
\begin{align*}
0 \to R\oplus P' \xrightarrow{1 \oplus f'} R \oplus Q' 
  \xrightarrow{(0,g')} A \to 0
\end{align*}
with some object $R\in \cal{P}$ and some minimal projective resolution
\begin{align*}
0 \to P' \xrightarrow{f'} Q' \xrightarrow{g'} A \to 0.
\end{align*}
\end{fct}

\begin{dfn}[{\cite[\S4.2]{B}}]\label{dfn:C_A}
Given an object $A\in \cal{A}$, 
take a minimal projective resolution 
\begin{align}
\label{eq:sec1:mpr}
0 \to P_A \xrightarrow{f_A} Q_A \xrightarrow{g} A \to 0,
\end{align}
We define a $\bb{Z}_2$-graded complex
\begin{align*}
\bu{{C_A}} 
 := \Bigl(\cpx{P_A}{f_A}{0}{Q_A}\Bigr) \in \cal{C}(\cal{P}).
\end{align*}
\end{dfn}

\begin{rmk}
The complex $\bu{{C_A}}$ is denoted as $C_A$ in \cite{B}.
\end{rmk}

By Fact~\ref{fct:B:lem4.1}, 
arbitrary two minimal projective resolutions of $A$ are isomorphic,
so the complex $\bu{{C_A}}$ is well-defined up to isomorphism.

\begin{fct}[{\cite[Lemma\,4.2]{B}}]
\label{fct:B:lem4.2}
Every object $M_\blob\in C(\cal{P})$ has a direct sum decomposition
\begin{align*}
\bu{M} = \bu{{C_A}} \oplus \bu{{C_B}}^* 
          \oplus \bu{{K_P}} \oplus \bu{{K_Q}}^*.
\end{align*}
Moreover, 
the objects $A,B \in \cal{A}$ and $P,Q \in \cal{P}$ 
are unique up to isomorphism.
\end{fct}

\subsubsection{Triangular decomposition}

\begin{dfn}[{\cite[\S\S4.3--4.4]{B}}]
\label{dfn:sec1:EF}
Given an object $A \in \cal{A}$, 
we define elements  $\bu{{E_A}},\bu{{F_A}} \in \cal{DH}(\cal{A})$ by
\begin{align*}
E_A :=  t^{\<\cl{P},\cl{A}\>} \cdot K_{-\cl{P}}*[\bu{{C_A}}],\quad
F_A :=  E_A^*.
\end{align*}
Here we used a projective decomposition \eqref{eq:sec1:mpr} of $A$
and the associated complex $\bu{{C_A}}$ in Definition~\ref{dfn:C_A}.
\end{dfn}

\begin{fct}[{\cite[Lemmas~4.6, 4.7]{B}}]
\begin{enumerate}
\item[(1)]
There is an embedding of algebras 
\begin{align*}
I_+^{\e}: \wtld{\cal{H}}(\cal{A}) \longinto \cal{DH}(\cal{A})
\end{align*}
defined by
$$
 [A] \longmapsto E_A \ (A\in\Iso(\cal{A})), \quad
 K_{\alpha} \longmapsto K_\alpha \ (\alpha\in K(\cal{A})).
$$
By composing $I_+^{\e}$ and the involution $*$, 
we also have an embedding 
\begin{align*}
I_{-}^{\e}: \wtld{\cal{H}}(\cal{A}) \longinto \cal{DH}(\cal{A})
\end{align*}
defined by 
$$
 [A] \longmapsto F_A \ (A \in \Iso(\cal{A})), \quad
 K_{\alpha} \longmapsto K_\alpha^* \ (\alpha \in K(\cal{A})).
$$

\item[(2)]
The multiplication map 
$$
 m: a \otimes b \longmapsto I_{+}^{\e}(a) * I_{-}^{\e}(b)
$$
defines  an isomorphism of vector spaces 
\begin{align*}
m: \wtld{\cal{H}}(\cal{A}) \otimes_{\bb{C}} \wtld{\cal{H}}(\cal{A}) 
     \xrightarrow{\ \sim\ } \cal{DH}(\cal{A}).
\end{align*}
\end{enumerate}
\end{fct}

As a corollary, we have 
\begin{cor}
\label{cor:sec1:basis}
$\cal{DH}(\cal{A})$ has a basis consisting of elements 
\begin{align*}
E_A * K_\alpha * K_{\beta}^* * F_B,\quad
A,B \in \Iso(\cal{A}),\ \alpha,\beta \in K(\cal{A}).
\end{align*}
\end{cor}


\subsection{Main theorem}

Now we can state our main theorem.

\begin{thm}
\label{thm:main}
Assume that the abelian category $\cal{A}$ satisfies 
the conditions (a)--(e).
Then the algebra $\cal{DH}(\cal{A})$ 
is isomorphic to the Drinfeld double of 
the bialgebra $\wtld{\cal{H}}(\cal{A})$.
\end{thm}

\section{Proof of the main theorem}
\label{sec:proof}

Because of the description 
of the basis of $\cal{DH}(\cal{A})$ 
(Corollary \ref{cor:sec1:basis})
and the definition of Drinfeld double
(Fact \ref{fct:sec1:dd}),
the proof of Theorem \ref{thm:main} 
is reduced to check the equation \eqref{eq:sec1:dd} 
for the elements 
consisting of the basis of $\wtld{\cal{H}}(\cal{A})$.
\\

Let us write the equation \eqref{eq:sec1:dd} in the present situation:
\begin{equation}
\label{eq:sec2:dd}
  \sum (a_{(2)},b_{(2)})_H \cdot 
   I_{+}^{\e}(a_{(1)}) * I_{-}^{\e}(b_{(1)})
 \stackrel{?}{=}
  \sum (a_{(1)},b_{(1)})_H \cdot 
   I_{-}^{\e}(b_{(2)}) * I_{+}^{\e}(a_{(2)}).
\end{equation}
What we must do is to check the equation 
for the cases 
\begin{equation}\label{eq:sec2:cases}
 (1)\ \ (a,b)=(K_{\alpha},K_{\beta}), \quad
 (2)\ \ (a,b)=([A],K_{\beta}), \quad
 (2')\ \ (a,b)=(K_{\alpha},[B]), \quad
 (3)\ \ (a,b)=([A],[B])
\end{equation}
with arbitrary $\alpha,\beta \in K(\cal{A})$ and 
$A, B \in \Iso(\cal{A})$.
\\

\subsection{Case (1)}

It is easy to check the equation \eqref{eq:sec2:dd}
for the case (1) in \eqref{eq:sec2:cases}.
Since $\Delta(K_{\alpha})=K_{\alpha} \otimes K_{\alpha}$,
the equation in this case becomes
$$
 (K_{\alpha},K_{\beta})_H  \cdot K_{\alpha} * K_{\beta}^*
  \stackrel{?}{=}
 (K_{\alpha},K_{\beta})_H  \cdot K_{\beta}^* * K_{\alpha},
$$ 
which is valid by Corollary \ref{cor:sec1:kkc} (2).

\subsection{Cases (2) and (2')}

The cases (2) and (2') in \eqref{eq:sec2:cases} are trivial.
In fact,
for the case (2),
we may write
$$
 \Delta([A])=\sum_{A_1,A_2} g_{A_1,A_2}^A  \cdot [A_1] \otimes [A_2]
$$
with some $g_{A_1,A_2}^A \in \bb{C}$ 
and $\Delta(K_\beta) = K_\beta \otimes K_\beta$.
Then \eqref{eq:sec2:dd} reads
$$
 \sum_{A_1,A_2} 
  ([A_2],K_{\beta})_H g_{A_1,A_2}^A  \cdot 
  E_{A_1} * K_{\beta}
 \stackrel{?}{=}
 \sum_{A_1,A_2} 
  ([A_1],K_{\beta})_H g_{A_1,A_2}^A  \cdot 
  K_{\beta} * E_{A_2}. 
$$
But recalling the Hopf pairing 
\begin{equation}
\label{eq:sec2:hp}
 ([A] K_\alpha,[B] K_\beta)_H 
  = \dfrac{\delta_{A,B}}{|\Aut_{\cal{A}}(A)|}t^{(\alpha,\beta)},
\end{equation}
in Definition/Fact \ref{df:sec1:hp},
the equation reads
$$
 g_{A,0}^A  \cdot E_{A} * K_{0}
 \stackrel{?}{=}
 g_{0,A}^A  \cdot K_{0} * E_{A}. 
$$
this equation trivially holds 
by Corollary \ref{cor:sec1:kkc}(1) and $g_{A,0}^A=g_{0,A}^A$.
The case (2') is similar.

\subsection{Case (3)}

The case (3) is non-trivial.
Let us write
\begin{align*}
 \Delta([A])&=
  \sum_{A_1,A_2} t^{\<A_1,A_2\>} 
  \dfrac{|\Ext_{\cal{A}}(A_1,A_2)_A|}{|\Aut_{\cal{A}}(A)|}  \cdot 
   [A_1] \otimes [A_2],
\\
 \Delta([B])&=
  \sum_{B_1,B_2} t^{\<B_1,B_2\>} 
  \dfrac{|\Ext_{\cal{A}}(B_1,B_2)_B|}{|\Aut_{\cal{A}}(B)|}  \cdot 
   [B_1] \otimes [B_2].
\end{align*}
Then the equation \eqref{eq:sec2:dd} reads
\begin{equation}
\label{eq:sec2:3:dd}
\begin{split}
&\sum_{A_1,A_2,B_1,B_2}
   t^{\<A_1,A_2\>} 
   \dfrac{|\Ext_{\cal{A}}(A_1,A_2)_A|}{|\Aut_{\cal{A}}(A)|}  
   t^{\<B_1,B_2\>} 
   \dfrac{|\Ext_{\cal{A}}(B_1,B_2)_B|}{|\Aut_{\cal{A}}(B)|}  
   (A_2,B_2)_H  \cdot 
   E_{A_1} * F_{B_1}
\\
\stackrel{?}{=}
&\sum_{A_1,A_2,B_1,B_2}
   t^{\<A_1,A_2\>} 
   \dfrac{|\Ext_{\cal{A}}(A_1,A_2)_A|}{|\Aut_{\cal{A}}(A)|}  
   t^{\<B_1,B_2\>} 
   \dfrac{|\Ext_{\cal{A}}(B_1,B_2)_B|}{|\Aut_{\cal{A}}(B)|}  
   (A_1,B_1)_H  \cdot 
   F_{B_2} * E_{A_2}.
\end{split}
\end{equation}

By the Hopf pairing \eqref{eq:sec2:hp},
the left hand side of \eqref{eq:sec2:3:dd} becomes
\begin{align*}
 \text{LHS of } \eqref{eq:sec2:3:dd}
&=\sum_{A_1,A_2,B_1,B_2}
   t^{\<A_1,A_2\>+\<B_1,B_2\>} 
   \dfrac{|\Ext_{\cal{A}}(A_1,A_2)_A|}{|\Aut_{\cal{A}}(A)|}  
   \dfrac{|\Ext_{\cal{A}}(B_1,B_2)_B|}{|\Aut_{\cal{A}}(B)|}  
   \dfrac{\delta_{A_2,B_2}}{|\Aut_{\cal{A}}(A_2)|}  \cdot 
   E_{A_1} * F_{B_1}
\\
&=\sum_{A_1,A_2,B_1}
   t^{\<A_1+B_1,A_2\>} 
   \dfrac{|\Ext_{\cal{A}}(A_1,A_2)_A|}{|\Aut_{\cal{A}}(A)|}  
   \dfrac{|\Ext_{\cal{A}}(B_1,A_2)_B|}{|\Aut_{\cal{A}}(B)|}  
   \dfrac{1}{|\Aut_{\cal{A}}(A_2)|}  \cdot 
   E_{A_1} * F_{B_1}.
\end{align*}
By Definition \ref{dfn:sec1:EF}, it becomes
\begin{align*}
&=\sum_{A_1,A_2,B_1}
   t^{\<A_1+B_1,A_2\>+\<P_{A_1},A_1\>+\<P_{B_1},B_1\>} 
   \dfrac{|\Ext_{\cal{A}}(A_1,A_2)_A|}{|\Aut_{\cal{A}}(A)|}  
   \dfrac{|\Ext_{\cal{A}}(B_1,A_2)_B|}{|\Aut_{\cal{A}}(B)|}  
   \dfrac{1}{|\Aut_{\cal{A}}(A_2)|}  
\\
&\phantom{=\sum_{A_1,A_2,B_1}}
   \cdot 
     K_{-\cl{P}_{A_1}}   * [\bu{{C_{A_1}}}] 
   * K_{-\cl{P}_{B_1}}^* * [\bu{{C_{B_1}}}^*].
\end{align*}
By Corollary \ref{cor:sec1:kkc} (1), it becomes
\begin{equation*}
\begin{split}
&=\sum_{A_1,A_2,B_1}
   t^{\<A_1+B_1,A_2\>+\<P_{A_1},A_1\>+\<P_{B_1},B_1\>
      -(P_{B_1},\cl{C}_{A_1 \bullet})} 
   \dfrac{|\Ext_{\cal{A}}(A_1,A_2)_A| \cdot |\Ext_{\cal{A}}(B_1,A_2)_B|}
         {|\Aut_{\cal{A}}(A)| \cdot |\Aut_{\cal{A}}(B)| 
          \cdot |\Aut_{\cal{A}}(A_2)|}  
\\
&\phantom{=\sum_{A_1,A_2,B_1}}
   \cdot 
     K_{-\cl{P}_{A_1}} * K_{-\cl{P}_{B_1}}^*
   * [\bu{{C_{A_1}}}]  * [\bu{{C_{B_1}}}^*].
\end{split}
\end{equation*}
Since $\cl{C}_{A_1 \bullet}=\cl{P}_{A_1}-\cl{Q}_{A_1}=\cl{A}$,
it becomes
\begin{equation}
\label{eq:sec2:3:LHS}
\begin{split}
&=\sum_{A_1,A_2,B_1}
   t^{\<A_1+B_1,A_2\>+\<P_{A_1},A_1\>+\<P_{B_1},B_1\>-(P_{B_1},A)} 
   \dfrac{|\Ext_{\cal{A}}(A_1,A_2)_A| \cdot |\Ext_{\cal{A}}(B_1,A_2)_B|}
         {|\Aut_{\cal{A}}(A)| \cdot |\Aut_{\cal{A}}(B)| 
          \cdot |\Aut_{\cal{A}}(A_2)|}  
\\
&\phantom{=\sum_{A_1,A_2,B_1}}
   \cdot 
     K_{-\cl{P}_{A_1}} * K_{-\cl{P}_{B_1}}^*
   * [\bu{{C_{A_1}}}]  * [\bu{{C_{B_1}}}^*].
\end{split}
\end{equation}

Now by the definition of multiplication 
\eqref{eq:sec1:diam} and \eqref{eq:sec1:*} in $\cal{DH}(\cal{A})$,
one may write
\begin{align}
\label{eq:sec2:3:LHS:ext}
[\bu{{C_{A_1}}}]  * [\bu{{C_{B_1}}}^*]
=\sum_{\bu{M} \in \Iso(\cal{C}(\cal{P}))} 
 t^{\<P_{A_1},Q_{B_1}\>+\<Q_{A_1},P_{B_1}\>}
 \dfrac{|\Ext_{\cal{C}(\cal{A})}(\bu{{C_{A_1}}},\bu{{C_{B_1}}}^*)_{\bu{M}}|}
       {|\Hom_{\cal{C}(\cal{A})}(\bu{{C_{A_1}}},\bu{{C_{B_1}}}^*)|}
 \cdot [\bu{M}].
\end{align}
Here the complex $\bu{M}$ sits in the commutative diagram
\begin{align*}
  \xymatrix{
      Q_{B_1} \ar@<.5ex>[r]^{0} \ar[d]
    &         \ar@<.5ex>[l]^{-f_{B_1}}  \ar[d] P_{B_1}
   \\
      M_1 \ar@<.5ex>[r]^{} \ar[d] 
    &     \ar@<.5ex>[l]^{} \ar[d] M_0
   \\
      P_{A_1} \ar@<.5ex>[r]^{f_{A_1}} 
    &     \ar@<.5ex>[l]^{0}  Q_{A_1}
  }
\end{align*}
where both columns give short exact sequences in $\cal{A}$.
Then, since $P_{A_1}$ and $Q_{A_1}$ are projective,
we have $M_1 \cong Q_{B_1} \oplus P_{A_1}$ and 
$M_0 \cong P_{B_1} \oplus Q_{A_1}$.
\begin{align}
\label{diag:sec2:M}
  \xymatrix{
      Q_{B_1} \ar@<.5ex>[r]^{0} \ar[d]
    &         \ar@<.5ex>[l]^{-f_{B_1}}  \ar[d] P_{B_1}
   \\
      Q_{B_1} \oplus P_{A_1} 
       \ar@<.5ex>[r]^{f_1}
       \ar[d] 
    &     \ar@<.5ex>[l]^{f_0} \ar[d] P_{B_1} \oplus Q_{A_1}
   \\
      P_{A_1} \ar@<.5ex>[r]^{f_{A_1}} 
    &     \ar@<.5ex>[l]^{0}  Q_{A_1}
  }
\end{align}
Now the commutativity of the diagram restricts 
the morphisms $M_1 \to M_0$ and $M_0 \to M_1$ of the following 
types:
\begin{align*}
f_1=\begin{pmatrix}0&s_{1}\\0&f_{A_1}\end{pmatrix},\quad
f_0=\begin{pmatrix}-f_{B_1}&s_{0}\\0&0\end{pmatrix}
\end{align*}
with $s_{1}: P_{A_1} \to P_{B_1}$ and  $s_{0}: Q_{A_1} \to Q_{B_1}$
(see also the argument in the proof of \cite[Lemma 3.3]{B}).
\\

Note also that 
\begin{equation*}
 \Hom_{\cal{C}(\cal{A})}(\bu{{C_{A_1}}},\bu{{C_{B_1}}}^*)
 \cong
 \Hom_{\cal{A}}(P_{A_1},Q_{B_1}),
\end{equation*}
which can be easily check.
Then the term 
$|\Hom_{\cal{C}(\cal{A})}(\bu{{C_{A_1}}},\bu{{C_{B_1}}}^*)|$ 
at the denominator of \eqref{eq:sec2:3:LHS:ext}
is equal to
\begin{equation}
\label{eq:sec2:3:LHS:ext:den}
 |\Hom_{\cal{C}(\cal{A})}(\bu{{C_{A_1}}},\bu{{C_{B_1}}}^*)|
=|\Hom_{\cal{A}}(P_{A_1},Q_{B_1})|
=t^{2\<P_{A_1},Q_{B_1}\>}
\end{equation}
At the last equation,
we used the hereditary property of $\cal{A}$ 
and the projective property of $\cal{P}_{A_1}$ and $\cal{Q}_{B_1}$.
\\

Combining \eqref{eq:sec2:3:LHS}, \eqref{eq:sec2:3:LHS:ext}
and \eqref{eq:sec2:3:LHS:ext:den},
the light hand side of \eqref{eq:sec2:3:dd} becomes
\begin{equation}
\label{eq:sec2:3:LHS:1}
\begin{split}
 \text{LHS of } \eqref{eq:sec2:3:dd}
&=\sum_{A_1,A_2,B_1,\bu{M}}
   t^{\<A_1+B_1,A_2\>+\<P_{A_1},A_1\>+\<P_{B_1},B_1\>
      -(P_{B_1},P_{A_1})+(P_{B_1},Q_{A_1})
      -\<P_{A_1},Q_{B_1}\>+\<Q_{A_1},P_{B_1}\>} 
\\
&\phantom{=\sum_{A_1,A_2,B_1}}
   \cdot 
   \dfrac{|\Ext_{\cal{A}}(A_1,A_2)_A| \cdot |\Ext_{\cal{A}}(B_1,A_2)_B|
          \cdot 
          |\Ext_{\cal{C}(\cal{A})}(\bu{{C_{A_1}}},\bu{{C_{B_1}}}^*)_{\bu{M}}|}
         {|\Aut_{\cal{A}}(A)| \cdot |\Aut_{\cal{A}}(B)| 
          \cdot |\Aut_{\cal{A}}(A_2)|}  
\\
&\phantom{=\sum_{A_1,A_2,B_1}}
   \cdot 
     K_{-\cl{P}_{A_1}} * K_{-\cl{P}_{B_1}}^*
   * [\bu{M}]
\end{split}
\end{equation}


A similar argument yields
\begin{equation}
\label{eq:sec2:3:RHS:1}
\begin{split}
 \text{RHS of } \eqref{eq:sec2:3:dd}
&=\sum_{\wtld{A}_1,\wtld{A}_2,\wtld{B}_2,\bu{N}}
   t^{\<\wtld{A}_1,\wtld{A}_2+\wtld{B}_2\>
      +\<P_{\wtld{B}_2},\wtld{B}_2\>+\<P_{\wtld{A}_2},\wtld{A}_2\>
      -(P_{\wtld{A}_2},\wtld{B}_2)
      +\<Q_{\wtld{B}_2},P_{\wtld{A}_2}\>-\<P_{\wtld{B}_2},Q_{\wtld{A}_2}\>} 
\\
&\phantom{=\sum_{\wtld{A}_1,\wtld{A}_2,\wtld{B}_1}}
   \cdot 
   \dfrac{|\Ext_{\cal{A}}(\wtld{A}_1,\wtld{A}_2)_A| 
          \cdot |\Ext_{\cal{A}}(\wtld{A}_1,\wtld{B}_2)_B|
          \cdot
          |\Ext_{\cal{C}(\cal{A})}
           (\bu{{C_{\wtld{B}_2}}}^*,\bu{{C_{\wtld{A}_2}}})_{\bu{N}}|}
         {|\Aut_{\cal{A}}(A)| \cdot |\Aut_{\cal{A}}(B)| 
          \cdot |\Aut_{\cal{A}}(A_1)|}  
\\
&\phantom{=\sum_{\wtld{A}_1,\wtld{A}_2,\wtld{B}_1}}
   \cdot 
     K_{-\cl{P}_{\wtld{A}_2}} * K_{-\cl{P}_{\wtld{B}_2}}^*
   * [\bu{N}]
\end{split}
\end{equation}
Here the complex $\bu{N}$ is of the next form:
\begin{align}
\label{diag:sec2:N}
  \xymatrix{
      P_{\wtld{A}_2} \ar@<.5ex>[r]^{f_{\wtld{A}_2}} \ar[d]
    &         \ar@<.5ex>[l]^{0}  \ar[d] Q_{\wtld{A}_2}
   \\
      P_{\wtld{A}_2} \oplus Q_{\wtld{B}_2} 
       \ar@<.5ex>[r]^{f_1'}
       \ar[d] 
    &  \ar@<.5ex>[l]^{f_0'} \ar[d] 
      Q_{\wtld{A}_2} \oplus P_{\wtld{B}_2}
   \\
      Q_{\wtld{B}_2} \ar@<.5ex>[r]^{0} 
    &     \ar@<.5ex>[l]^{-f_{\wtld{B}_2}}  P_{\wtld{B}_2}
  }
\end{align}
where
\begin{align*}
f_1'=\begin{pmatrix}f_{\wtld{A}_2}&s_1'\\0&0\end{pmatrix},\quad
f_0'=\begin{pmatrix}0&s_0'\\0&-f_{\wtld{B}_2}\end{pmatrix}
\end{align*}
with $s_{1}': Q_{\wtld{B}_2} \to Q_{\wtld{A}_2}$ 
and  $s_{0}': P_{\wtld{B}_2} \to P_{\wtld{A}_2}$.
\\

Now we must compare the coefficients of 
the same term $[\bu{M}]=[\bu{N}]$ 
in \eqref{eq:sec2:3:LHS:1} and \eqref{eq:sec2:3:RHS:1}.
A quick observation yields that we must have the correspondences
\begin{equation}
\label{eq:sec2:3:corres1}
 Q_{\wtld{A}_2} = Q_{B_1}, \quad
 Q_{\wtld{B}_2} = P_{B_1}, \quad
 P_{\wtld{B}_2} = P_{A_1}, \quad
 P_{\wtld{A}_2} = Q_{A_1}
\end{equation}
and
\begin{equation}
\label{eq:sec2:3:corres2}
 f_{\wtld{B}_2} = -s_1, \quad 
 f_{\wtld{A}_2} = s_0, \quad 
 s_1' = -f_{B_1}, \quad 
 s_0' = f_{A_1}. 
\end{equation}
Then from \eqref{diag:sec2:M} and \eqref{diag:sec2:N},
we have the next combined diagram
\begin{align*}
\xymatrix{
   P_{A_1} \ar@{^{(}_->}[r]^{f_{A_1}} 
           \ar@{^{(}_->}[d]_{-f_{\wtld{B}_2}=s_1}
 & Q_{A_1} \ar@{->>}[r]  \ar[d]^{f_{\wtld{A}_2}=s_0'}
 & A_1
 \\
   P_{B_1} \ar@{^{(}_->}[r]_{f_{B_1}} \ar@{->>}[d]
 & Q_{B_1} \ar@{->>}[r]  \ar@{->>}[d]
 & B_1
 \\
   \wtld{B}_1 
 & \wtld{A}_2
}
\end{align*}
where all the columns and rows are short exact. 
We also have the short exact sequences
\begin{align*}
\xymatrix{
   A_2 \ar@{^{(}_->}[r] \ar@{=}[d]
 & A \ar@[->>][r]
 & A_1
 & \wtld{A}_2 \ar@{^{(}_->}[r] 
 & A \ar@[->>][r]
 & \wtld{A}_1 \ar@{=}[d]
 \\
   B_2 \ar@{^{(}_->}[r] 
 & B \ar@[->>][r]
 & B_1
 & \wtld{B}_2 \ar@{^{(}_->}[r] 
 & B \ar@[->>][r]
 & \wtld{B}_1
}
\end{align*}
Considering the classes in the Grothendieck group $K(\cal{A})$,
we have the relations
\begin{align*}
 &\cl{\wtld{A}}_1=\cl{A}-\cl{\wtld{A}}_2 
  = \cl{A}-(\cl{Q}_{B_1}-\cl{Q}_{A_1})
\\
=&\cl{\wtld{B}}_1=\cl{B}-\cl{\wtld{B}}_2 
  = \cl{B}-(\cl{P}_{B_1}-\cl{P}_{A_1})
\end{align*}
and
\begin{align*}
 & \cl{A}_2=\cl{A}-\cl{A}_1 
  = \cl{A}-(\cl{Q}_{A_1}-\cl{P}_{A_1})
\\
=& \cl{B}_2=\cl{B}-\cl{B}_1 
  = \cl{B}-(\cl{Q}_{B_1}-\cl{P}_{B_1}).
\end{align*}
Then we have
\begin{align*}
 \cl{A}-\cl{B} 
 &=(\cl{Q}_{B_1}-\cl{Q}_{A_1})-(\cl{P}_{B_1}-\cl{P}_{A_1})
\\
 &=(\cl{Q}_{A_1}-\cl{P}_{A_1})-(\cl{Q}_{B_1}-\cl{P}_{B_1}),
\end{align*}
so that in $K(\cal{A})$ we have
\begin{equation}
\label{eq:sec2:3:cls}
 \cl{Q}_{B_1} + \cl{P}_{A_1} = \cl{Q}_{A_1} + \cl{P}_{B_1},\quad
 \cl{A}_1 = \cl{B}_1, \quad
 \cl{\wtld{A}}_2 = \cl{\wtld{B}}_2.
\end{equation}

Now we will finish the proof.
By the correspondences \eqref{eq:sec2:3:corres1} 
and \eqref{eq:sec2:3:corres2},
we have only to check that the coefficients 
\begin{equation}
\label{eq:sec2:3:LHS:2}
\<A_1+B_1,A_2\>+\<P_{A_1},A_1\>+\<P_{B_1},B_1\>
-(P_{B_1},A_1)
-\<P_{A_1},Q_{B_1}\>+\<Q_{A_1},P_{B_1}\>
\end{equation}
in \eqref{eq:sec2:3:LHS:1} and 
\begin{equation}
\label{eq:sec2:3:RHS:2}
\<\wtld{A}_1,\wtld{A}_2+\wtld{B}_2\>
+\<P_{\wtld{B}_2},\wtld{B}_2\>+\<P_{\wtld{A}_2},\wtld{A}_2\>
-(P_{\wtld{A}_2},\wtld{B}_2)
+\<Q_{\wtld{B}_2},P_{\wtld{A}_2}\>-\<P_{\wtld{B}_2},Q_{\wtld{A}_2}\>
\end{equation}
in \eqref{eq:sec2:3:RHS:1} coincide.
But using \eqref{eq:sec2:3:cls} we have
\begin{align*}
\eqref{eq:sec2:3:LHS:2}
&=\<A_1+B_1,A\>-\<A_1+B_1,A_1\>+\<P_{A_1},A_1\>+\<P_{B_1},B_1\>
\\
&\phantom{=<}
-\<P_{B_1},A_1\>-\<A_1,P_{B_1}\>
-\<P_{A_1},Q_{B_1}\>+\<Q_{A_1},P_{B_1}\>
\\
&=2\<A_1,A_2\>+\<P_{A_1},A_1\>-\<A_1,P_{B_1}\>
  -\<P_{A_1},Q_{B_1}\>+\<Q_{A_1},P_{B_1}\>
\\
&=2\<A_1,A_2\>+\<P_{A_1},\cl{A}_1-\cl{Q}_{B_1}\>
 +\<\cl{Q}_{A_1}-\cl{A}_1,P_{B_1}\>
\\
&=2\<A_1,A_2\>+\<P_{A_1},-\cl{P}_{B_1}\>
 +\<\cl{P}_{A_1},P_{B_1}\>
=2\<A_1,A_2\>.
\end{align*}
A similar calculation gives
\begin{align*}
\eqref{eq:sec2:3:RHS:2}=2\<A_1,A_2\>.
\end{align*}
Thus the proof is completed.


\section{Concluding remarks}
\label{sec:cr}

As mentioned in \cite[\S1.4]{B},
the work of Cramer \cite{C} and Theorem \ref{thm:main}
yield  the following:

\begin{cor}
For the hereditary abelian cateogry,
the algebra $\cal{DH}(\cal{A})$ is functorial 
with respect to derived invariance.

Precisely speaking,
let $\cal{A}$ and $\cal{B}$ 
be two abelian categories satisfying conditions (a)--(e).
Assume that the bounded derived categories 
$D^b(\cal{A})$ and $D^b(\cal{B})$
of $\cal{A}$ and $\cal{B}$ are equivalent by the functor $\Phi$:
$$
 \Phi: D^b(\cal{A}) \xrightarrow{\ \sim \ } D^b(\cal{B}).
$$
Then one can construct an algebra isomorphism
$$
 \Phi^{\cal{DH}}: 
 \cal{DH}(\cal{A}) \xrightarrow{\ \sim\ } \cal{DH}(\cal{B}),
$$
and this construction is functorial:
$(\Phi_1 \circ \Phi_2)^{\cal{DH}} = \Phi_1^{\cal{DH}} \circ \Phi_2^{\cal{DH}}$
for all equivalences $\Phi_1, \Phi_2$.
\end{cor}

Now set $\cal{T}:=D^b(\cal{A})$ 
for some hereditary abelian category $\cal{A}$ satisfying (a)--(e).
We also set 
$$
 \cal{DH}(\cal{T}) := \cal{DH}(\cal{A}).
$$
This algebra depends only on the triangulated cateogry $\cal{T}$
by the above corollary.
Denote by $\Auteq(\cal{T})$
the group of autoequivalences of $\cal{T}$.
Then, setting 
$$
  \Aut^{\cal{DH}}(\cal{T}) := 
  \{\Phi^{\cal{DH}} \mid \Phi \in \Auteq(\cal{T}) \},
$$
we have an embedding of groups
$$
 \Aut^{\cal{DH}}(\cal{T}) \subset \Aut(\cal{DH}(T)).
$$
\\

Let us close this note by mentioning a non-trivial example.
For an elliptic curve $C$ defined over $\frk{k}$,
set
$$
 \cal{A} := \Coh(C),
$$
the abelian category of coherent sheaves on $C$.
This category satisfies the conditions (a)--(e).
Set $\cal{T}:=D^b(\cal{A})$ as above.

Then by the theory of Fourier-Mukai transforms,
we have a short exact sequence 
$$
 \xymatrix{
   0 \ar[r] 
 & \bb{Z} \oplus (C \times \wh{C}) \ar[r] 
 & \Auteq(\cal{T}) \ar[r] 
 & \SL_2(\bb{Z}) \ar[r]
 & 0
 }
$$
of groups.
Here $\bb{Z} \oplus (C \times \wh{C})$
corresponds to the subgroup of $\Auteq(\cal{T})$
generated by
the shifts $[n]$ of complexes,
the pushforward $t_{a*}$ by translations on $C$ with $a \in C$,
and tensor products $L \otimes (-)$ 
with $L \in \wh{C}:=\Pic^0(C)$.
The cokernel part $\SL_2(\bb{Z})$ 
consists of (non-trivial) Fourier-Mukai transforms 
$\Phi_{\cal{E}}
 := \bb{R}p_{2*}(\cal{E} \stackrel{\bb{L}}{\otimes} p_1^*(-))$
with $\cal{E} \in D^b(C \otimes \wh{C})$.
The generators 
$S:=\begin{pmatrix}0&-1\\1&0\end{pmatrix}$ and 
$T:=\begin{pmatrix}1&-1\\0&1\end{pmatrix}$ 
correspond to the Fourier-Mukai transforms 
$\Phi_{\cal{E}_0}$ and $L \otimes (-)$,
where $\cal{E}_0$ is the Poincar\'{e} bundle 
on $C \otimes \wh{C}$,
and $L \in \Pic(C)$ is a degree one line bundle.
(See \cite[\S9.5]{H} for the detailed explanation.)

Now one can see that 
the operation 
$$
 DH: \Phi \longmapsto \Phi^{\cal{DH}}
$$ 
yields
$$
 \xymatrix{
   0 \ar[r] 
 & \bb{Z} \oplus (C \times \wh{C}) \ar[r] \ar[d]^{DH}
 & \Auteq(\cal{T}) \ar[r] \ar[d]^{DH}
 & \SL_2(\bb{Z})   \ar[r] \ar[d]^{DH}
 & 0
 \\ 
   0 \ar[r] 
 & \bb{Z}/2\bb{Z} \ar[r]
 & \Aut^{\cal{DH}}(\cal{T}) \ar[r] 
 & \SL_2(\bb{Z}) \ar[r]
 & 0
 }
$$
Here $\bb{Z}/2\bb{Z}$ corresponds to the involution $*$
of the algebra $\cal{DH}(\cal{T})$.
Thus $\SL_2(\bb{Z})$ acts on $\cal{DH}(\cal{T})$.
This action is essentially the same as  
the $\SL_2(\bb{Z})$-automorphisms of the algebra 
appearing in \cite{BS}.
The same algebra appeared in the works 
\cite{FF}, \cite{FH}, \cite{FT} and \cite{SV1}.
The $\SL_2(\bb{Z})$-automorphisms 
(precisely speaking, the counterpart in the degenerate algebra) 
play an important role 
in the argument of \cite{SV2} 
in the context of the so-called AGT relation/conjecture.


\end{document}